# THE LARGEST COMPONENT IN A SUBCRITICAL RANDOM GRAPH WITH A POWER LAW DEGREE DISTRIBUTION

BY SVANTE JANSON

*Uppsala University*

It is shown that in a subcritical random graph with given vertex degrees satisfying a power law degree distribution with exponent $\gamma > 3$, the largest component is of order $n^{1/(\gamma-1)}$. More precisely, the order of the largest component is approximatively given by a simple constant times the largest vertex degree. These results are extended to several other random graph models with power law degree distributions. This proves a conjecture by Durrett.

**1. Introduction.** Random graphs where the asymptotic distribution of the vertex degrees has a power law tail have been the focus of much interest during the last decade or so, including both theoretical studies and various applications; see, for example, the books by Durrett [7] and (from a physical point of view) Dorogovtsev and Mendes [6].

To be precise, we will in this paper consider uniformly distributed random graphs with a given degree sequence, defined as follows. (See Section 4 for extensions to some other random graph models.) Let $n \in \mathbb{N}$ and let $(d_i)_1^n$ be a sequence of nonnegative integers. We let $G(n,(d_i)_1^n)$ be a random graph with degree sequence $(d_i)_1^n$, uniformly chosen among all possibilities (tacitly assuming that there is any such graph at all; in particular, $\sum_i d_i$ has to be even).

We consider asymptotics as $n \to \infty$, and thus we assume throughout the paper that we are given a sequence $(d_i)_1^n$ for each $n \in \mathbb{N}$ (or at least for some sequence $n \to \infty$); $d_i$ may depend on $n$ but for notational simplicity we will not show the dependency on $n$ explicitly. (Similarly, we often omit the index $n$ on other parameters and variables below.)

We let the random variable $D = D_n$ be the degree of a randomly chosen vertex in $G(n,(d_i)_1^n)$. Thus $\mathbb{P}(D = k) = n_k/n$, where $n_k := |\{i : d_i = k\}|$ is the number of vertices of degree $k$ in $G(n,(d_i)_1^n)$.

We say that $G(n,(d_i)_1^n)$ has a *power law degree distribution with exponent* $\gamma$ if $\mathbb{P}(D_n = k) \to p_k$ as $n \to \infty$ for every $k \geq 0$ and some probability distribution $(p_k)_0^\infty$, and this asymptotic degree distribution satisfies $p_k \sim ck^{-\gamma}$ as $k \to \infty$ for some $c > 0$. [This is thus a doubly asymptotic notion, letting first $n \to \infty$ and then $k \to \infty$. Actually, the theorems below will not assume this; we will use the weaker inequality (1.3) instead, but, on the other hand, we have to impose uniformity in $n$. Nevertheless, the reader ought to think mainly of the case of a power law degree distribution.] Note that necessarily $\gamma > 1$, since otherwise $(p_k)$ would not be summable.

The existence of a giant component (i.e., a component of order $n$) in $G(n,(d_i)_1^n)$ has been studied by Molloy and Reed [15, 16]; see also Janson and Łuczak [11]. To state their result we let $\mathcal{C}_1$ be the largest component of $G(n,(d_i)_1^n)$ and $|\mathcal{C}_1|$ its order. We further introduce the two parameters

$$(1.1) \qquad \mu = \mu_n := \mathbb{E} D_n = \frac{1}{n} \sum_{k=0}^\infty k n_k = \frac{1}{n} \sum_{i=1}^n d_i,$$

$$(1.2) \qquad \nu = \nu_n := \frac{1}{\mu n} \sum_{k=0}^\infty k(k-1) n_k = \frac{1}{n\mu} \sum_{i=1}^n d_i(d_i - 1);$$

$\mu$ is thus the average degree and $\nu$ is the mean number of offspring in the usual branching process approximation of the local structure starting at any fixed vertex; see, for example, [7]. Loosely speaking, the condition by Molloy and Reed [15] for existence of a giant component is $\nu > 1$. More precisely, if we assume as above that $n_k/n \to p_k$, $k \geq 0$, then, under suitable conditions ensuring uniform summability, $\mu_n \to \mu_\infty := \sum_{k=0}^\infty k p_k$ and $\nu_n \to \nu_\infty := \mu_\infty^{-1} \sum_{k=0}^\infty k(k-1) p_k$, and (under weak additional technical assumptions; see [11, 15]), if $\nu_\infty > 1$, then there is a constant $\rho > 0$ such that $|\mathcal{C}_1| = \rho n + o_p(n)$, while if $\nu_\infty \leq 1$, then $|\mathcal{C}_1| = o_p(n)$. (We use $o_p$ and $O_p$ in the standard way; see, e.g., Janson, Łuczak and Ruciński [12]. For example, $|\mathcal{C}_1| = o_p(n)$ means that $\mathbb{P}(|\mathcal{C}_1| > \varepsilon n) \to 0$ as $n \to \infty$ for every $\varepsilon > 0$; equivalently, $|\mathcal{C}_1|/n \xrightarrow{\mathrm{P}} 0$.)

Returning to power law degree distributions, we see that if $1 < \gamma \leq 2$, then $\mu_\infty = \infty$, and if $2 < \gamma \leq 3$, then $\mu_\infty < \infty$ and $\nu_\infty = \infty$; in both cases we expect a giant component. (We will not study these cases further.) In the sequel, we will assume $\gamma > 3$, and thus $\nu_\infty < \infty$. (Note that in this range, the condition $\nu_\infty > 1$ is determined mainly by the values of $p_k$ for small $k$, such as $k = 1$, and is essentially independent of the behavior of $p_k$ for large $k$ and in particular of the value of the exponent $\gamma$. Relations between $\gamma$ and the existence of a giant component that have been reported in the literature are artifacts due to the assumption of specific forms of $p_k$ also for small $k$.)

We will in this paper consider the subcritical case $\nu_\infty < 1$, when $|\mathcal{C}_1|$ is of smaller order than $n$, and try to estimate its order more precisely. Durrett [7] observed that while it is well known that for a subcritical Erdős–Rényi random graph $|\mathcal{C}_1| = O_p(\log n)$ (see, e.g., [1], Chapter V.2), this fails for graphs with power law degree distributions. In fact, typically the largest vertex degree is of order $n^{1/(\gamma-1)}$, and obviously the largest component is at least as big. Durrett ([7], Conjecture 3.3.1), conjectured (for a slightly different random graph model; see Section 4.1) that this is the correct order of the largest component. In the present paper, we prove this conjecture. A related result has recently been shown by Pittel [20], who showed that under the stronger condition $\mathbb{P}(D_n = k) = O(k^{-\gamma})$, which implies that the maximum degree is $O(n^{1/\gamma})$, the conclusion $|\mathcal{C}_1| = O_p(n^{1/\gamma} \log n)$ holds.

All unspecified limits below are as $n \to \infty$. We say that an event holds whp (*with high probability*), if it holds with probability tending to 1 as $n \to \infty$. We let $C_1, C_2, \ldots$ denote unspecified constants, sometimes depending on other parameters but never on $n$.

THEOREM 1.1. *Consider the random graph $G(n, (d_i)_1^n)$, with $(d_i)_1^n = (d_i^{(n)})_1^n$ as above, and let $\mathcal{C}_1$ be its largest component. Suppose, using the notation above, that $\mu_n \to \mu_\infty > 0$ and $\nu_n \to \nu_\infty < 1$. Suppose further that*

$$\mathbb{P}(D_n \geq k) = O(k^{1-\gamma}), \tag{1.3}$$

*uniformly in $n$ and $k \geq 1$, for some $\gamma > 3$. Then there exists a constant $A$ such that $|\mathcal{C}_1| \leq An^{1/(\gamma-1)}$ whp.*

REMARK 1.2. It is not necessary to assume that a limit distribution $(p_k)$ exists, and the assumptions $\mu_n \to \mu_\infty > 0$ and $\nu_n \to \nu_\infty < 1$ may be replaced by $\liminf \mu_n > 0$ and $\limsup \nu_n < 1$. The constant $A$ depends only on $\liminf \mu_n$, $\limsup \nu_n$ and the constant $C_1$ implicit in (1.3).

Let $\Delta = \Delta_n := \max_i d_i$, the maximum degree of $G(n, (d_i)_1^n)$. Note that the assumption (1.3) implies

$$n^{-1} \leq \mathbb{P}(D \geq \Delta) \leq C_1 \Delta^{1-\gamma}$$

and thus

$$\Delta \leq C_2 n^{1/(\gamma-1)}. \tag{1.4}$$

We can sharpen Theorem 1.1 as follows.

THEOREM 1.3. *Under the assumptions of Theorem 1.1,*

$$|\mathcal{C}_1| = \frac{\Delta}{1-\nu} + o_p(n^{1/(\gamma-1)}). \tag{1.5}$$

As said above, typically we expect $\Delta = \Delta_n$ to be of the order $n^{1/(\gamma-1)}$, and then Theorem 1.3 shows that $|\mathcal{C}_1|$ is of this order also. The factor $1/(1-\nu)$ in (1.5) is exactly what the standard branching process approximation of the local structure suggests: If we start at a vertex $v^*$ with degree $\Delta$, it has $\Delta$ neighbors; each of them has on the average about $\nu$ further neighbors, so there are about $\nu\Delta$ vertices of distance 2 from $v^*$. Continuing, we expect to find about $\nu^{j-1}\Delta$ vertices of distance $j$ from $v^*$, and thus in total about $\Delta/(1-\nu)$ vertices in the component containing $v^*$; see further Section 3.

We extend Theorems 1.1 and 1.3 to some other random graph models in Section 4.

We have not attempted to find any sharper estimate of the error term in (1.5), and leave this as an open problem for the reader. This is especially interesting when $\Delta = o(n^{1/(\gamma-1)})$. For example, as said above, Pittel [20] studies a special case of the situation above with $\Delta = O(n^{1/\gamma})$, and shows that under his assumptions $|\mathcal{C}_1| = O_p(n^{1/\gamma} \log n)$.

REMARK 1.4. The proof in Section 3 is easily modified to show that if $\mathcal{C}_j$ is the $j$th largest component and $\Delta^{(j)}$ is the $j$th largest vertex degree, then

$$|\mathcal{C}_j| = \frac{\Delta^{(j)}}{1-\nu} + o_p(n^{1/(\gamma-1)})$$

for every fixed $j \geq 1$. Moreover, it follows that, for every fixed $J$, whp the $J$ vertices with largest degrees belong to different components, each of them being of order roughly $(1-\nu)^{-1}$ times the degree of that vertex. If these vertex degrees are sufficiently well separated, the largest component thus contains the vertex with largest degree, and so on (for at least $J-1$ components), but, of course, this can fail if two of the degrees are very close to each other.

REMARK 1.5. As is customary, we have concentrated on the number of vertices in the components in these results. The proofs show that the same results hold also if we consider the number of edges in the components. Consequently, the difference between the numbers of edges and vertices in the largest component is $o_p(n^{1/(\gamma-1)})$.

**2. Proof of Theorem 1.1.** It will be convenient to work with *multigraphs*, that is, to allow multiple edges and loops. More precisely, we shall use the following standard type of random multigraph: Let $n \in \mathbb{N}$ and let $(d_i)_1^n$ be a sequence of nonnegative integers such that $\sum_{i=1}^n d_i$ is even. We let $G^*(n, (d_i)_1^n)$ be the *random multigraph with given degree sequence* $(d_i)_1^n$, defined by the configuration model (see, e.g., Bollobás [1]): take a set of $d_i$ *half-edges* for each vertex $i$, and combine the half-edges into pairs by a uniformly random matching of the set of all half-edges (this pairing is called a *configuration*); each pair of half-edges is then joined to form an edge of $G^*(n,(d_i)_1^n)$.

We obtain $G(n,(d_i)_1^n)$ by conditioning the multigraph $G^*(n,(d_i)_1^n)$ on being a (simple) graph. It is well known that the conditions in Theorem 1.1 imply that

$$(2.1) \qquad \liminf_{n\to\infty} \mathbb{P}(G^*(n,(d_i)_1^n) \text{ is simple}) > 0,$$

and thus it suffices to prove the result for the random multigraph $G^*(n,(d_i)_1^n)$. (In fact, by Janson [9], a necessary and sufficient condition for (2.1) is $\mathbb{E}D_n^2 = O(\mathbb{E}D_n)$, or equivalently $\nu_n = O(1)$.)

We thus consider the random multigraph $G^*(n,(d_i)_1^n)$, and prove Theorem 1.1 for this case. Fix a vertex $v$ and explore the component $\mathcal{C}(v)$ containing $v$ by the standard exploration process. This process starts by declaring $v$ *used* and all half-edges at $v$ *active*. Then, proceed as follows for $i = 1, 2, \ldots$, as long as there is any active half-edge: Take an active half-edge, say $x_i$. (In this proof it does not matter which one we choose; we may use any deterministic or random rule. Later it will be convenient to use the breadth-first version where the active half-edges are processed in order of appearance.) Find its partner $y_i$ in the configuration [thus $x_i y_i$ is an edge in $G^*(n,(d_i)_1^n)$], and let $v_i$ be the other endpoint of the edge $x_i y_i$, that is, the endpoint of $y_i$. If the vertex $v_i$ is not already used, declare the remaining $d(v_i) - 1$ half-edges at $v_i$ active. Finally, declare $x_i$, $y_i$ and $v_i$ used. Repeat.

Let $S_i$ be the number of active half-edges after $i$ steps. Thus $S_0 = d(v)$ and $S_i = S_{i-1} + \xi_i - 1$, $i \geq 1$, where $\xi_i$ is the number of new half-edges found at step $i$. We have $\xi_i = d(v_i) - 1$ if $v_i$ is not already used; otherwise $\xi_i = 0$.

The process stops at $\tau := \min\{i : S_i = 0\}$; we then have found $\tau$ edges, so the component $\mathcal{C}(v)$ has $\tau$ edges and at most $\tau + 1$ vertices. For definiteness, we let $\xi_i = 0$ for $i > \tau$.

In particular, for any integer $M \geq 0$, if $|\mathcal{C}(v)| > M$, then $\tau \geq M$ so the process lives at least until $M$ and $0 \leq S_M = d(v) + \sum_{i=1}^M (\xi_i - 1)$ and thus

$$(2.2) \qquad d(v) + \sum_{i=1}^M \xi_i \geq M.$$

We reveal the random configuration only as we need it. Equivalently, we may construct the random configuration during the exploration, by choosing each partner $y_i$ uniformly among all half-edges except $x_i$ and the ones already paired. The numbers $\xi_i$ then become random variables. The random variables $\xi_i$ are dependent, but we can approximate them by independent random variables. More precisely, we bound them as follows.

Suppose that $i \leq \sqrt{n}$, say. (Our components will be much smaller.) If $\xi_i = k \geq 1$, then $d(v_i) = k + 1$, so there are at most $n_{k+1} = n\mathbb{P}(D = k+1)$ possible choices for $v_i$ (excluding the used ones) and thus at most $(k + $ $1)n_{k+1} = n(k+1)\mathbb{P}(D = k+1)$ choices of the half-edge $y_i$. When choosing $y_i$, we have revealed the pairings of $2(i-1)$ half-edges and chosen $x_i$, but $y_i$ is uniformly distributed over all $\sum_{j=1}^{n} d_j - 2(i-1) - 1 = n\mu - O(n^{1/2})$ remaining half-edges. Hence, given any history of the exploration process up to step $i$,

$$
(2.3) \quad \begin{aligned} \mathbb{P}(\xi_i = k) &\leq \frac{n(k+1)\mathbb{P}(D = k+1)}{n\mu - O(n^{1/2})} \\ &= \frac{(k+1)\mathbb{P}(D = k+1)}{\mu}(1 + O(n^{-1/2})). \end{aligned}
$$

Fix $\varepsilon > 0$ with $4\varepsilon \leq 1 - \nu_\infty$, and let $\nu' = \nu + \varepsilon$. We consider only $n$ that are so large that $\nu = \nu_n < 1 - 3\varepsilon$, and thus $\nu' < 1 - 2\varepsilon$. Let $X$ be an integer-valued random variable with the distribution given by

$$
(2.4) \quad \mathbb{P}(X \geq x) = \max\left(1, \frac{\nu'}{\nu} \sum_{k \geq x} \frac{(k+1)\mathbb{P}(D = k+1)}{\mu}\right), \quad x > 0.
$$

If $n$ is large enough (not depending on $i$ or $k$), the factor $1 + O(n^{-1/2})$ in (2.3) is less than $\nu'/\nu$, and thus by (2.3) and (2.4), given any previous history,

$$
\mathbb{P}(\xi_i \geq x) \leq \mathbb{P}(X \geq x) \quad \text{for every } x \geq 0.
$$

In other words, $\xi$ is stochastically dominated by $X$. In particular, this holds given $\xi_1, \ldots, \xi_{i-1}$, and it follows that the sequence $\xi_1, \xi_2, \ldots$ may be coupled with a sequence of independent random variables $X_i$ with $X_i \stackrel{\mathrm{d}}{=} X$, such that

$$
(2.5) \quad \xi_i \leq X_i, \quad 1 \leq i \leq \sqrt{n}.
$$

Furthermore, trivially,

$$
(2.6) \quad \xi_i \leq \Delta, \quad i \geq 1.
$$

The remainder of the proof will only use the properties (2.5) and (2.6) of $(\xi_i)$.

We note first that, by (2.4) and the assumptions of the theorem,

$$
(2.7) \quad \begin{aligned} \mathbb{E}X &= \sum_{m=1}^{\infty} \mathbb{P}(X \geq m) \leq \frac{\nu'}{\nu} \sum_{m=1}^{\infty} \sum_{k=m}^{\infty} \frac{(k+1)\mathbb{P}(D = k+1)}{\mu} \\ &= \frac{\nu'}{\nu} \sum_{k=1}^{\infty} k \frac{(k+1)\mathbb{P}(D = k+1)}{\mu} = \frac{\nu'}{\nu}\nu = \nu' \end{aligned}
$$

and, for all $x > 0$, with $m := \lceil x \rceil$,

$$\mathbb{P}(X \geq x) \leq \frac{\nu'}{\nu} \sum_{k=m}^{\infty} \frac{(k+1)\mathbb{P}(D = k+1)}{\mu}$$

$$\leq C_3 \sum_{k=m+1}^{\infty} k\mathbb{P}(D = k)$$

(2.8)
$$= C_3 \left( m\mathbb{P}(D > m) + \sum_{j=m+1}^{\infty} \mathbb{P}(D \geq j) \right)$$

$$\leq C_4 m^{2-\gamma} + C_4 \sum_{j=m+1}^{\infty} j^{1-\gamma} \leq C_5 m^{2-\gamma} \leq C_5 x^{2-\gamma}.$$

Let $M = An^{1/(\gamma-1)}$ and $M_1 = n^{1/(\gamma-1)-\delta}$ for some positive constants $A$ (large) and $\delta$ [small, and in particular less than $1/(\gamma-1)$] to be chosen later, and let $Y_i := X_i \mathbf{1}[X_i \leq M_1]$. Thus, by (2.7), $\mathbb{E}Y_i \leq \mathbb{E}X \leq \nu'$.

If $|\mathcal{C}(v)| > M$, then (2.2) holds and thus, using (2.6),

$$M \leq d(v) + \sum_{i=1}^{M} \xi_i \leq \Delta + \sum_{i=1}^{M} \xi_i \mathbf{1}[X_i \leq M_1] + \sum_{i=1}^{M} \xi_i \mathbf{1}[X_i > M_1]$$

(2.9)
$$\leq \sum_{i=1}^{M} Y_i + \Delta \left( 1 + \sum_{i=1}^{M} \mathbf{1}[X_i > M_1] \right)$$

$$\leq \nu' M + \sum_{i=1}^{M}(Y_i - \mathbb{E}Y_i) + \Delta \left( 1 + \sum_{i=1}^{M} \mathbf{1}[X_i > M_1] \right)$$

and thus, recalling $\nu' < 1 - 2\varepsilon$,

(2.10)
$$\mathbb{P}(|\mathcal{C}(v)| > M) \leq \mathbb{P}\left( \sum_{i=1}^{M}(Y_i - \mathbb{E}Y_i) > \varepsilon M \right)$$

$$+ \mathbb{P}\left( \sum_{i=1}^{M} \mathbf{1}[X_i > M_1] \geq \varepsilon \frac{M}{\Delta} - 1 \right).$$

Note that both sums on the right-hand side are sums of i.i.d. random variables, so we can, and shall, use simple standard estimates for them.

For the first sum in (2.10), fix a number $r \geq \gamma$ such that $r\delta > 2$. Then, by Rosenthal's inequality ([8], Theorem 3.9.1) (with constants depending on $r$),

$$\mathbb{E}\left| \sum_{i=1}^{M}(Y_i - \mathbb{E}Y_i) \right|^r \leq C_6 M^{r/2}(\mathbb{E}|Y_1 - \mathbb{E}Y_1|^2)^{r/2} + C_7 M \mathbb{E}|Y_1 - \mathbb{E}Y_1|^r$$

(2.11)
$$\leq C_6 M^{r/2}(\mathbb{E}Y_1^2)^{r/2} + C_8 M \mathbb{E}Y_1^r.$$

We estimate the moments of $Y_1$ and obtain using (2.8)

(2.12)
$$\begin{aligned}\mathbb{E}Y_1^2 &= \int_0^\infty 2x\mathbb{P}(Y_1 > x)\,dx \\ &= \int_0^{M_1} 2x\mathbb{P}(X > x)\,dx \\ &\leq 1 + C_5 \int_1^{M_1} 2x^{3-\gamma}\,dx \leq C_9 M_1\end{aligned}$$

(rather coarsely and because $\gamma > 3$), and trivially

(2.13)
$$\mathbb{E}Y_1^r \leq M_1^r.$$

Consequently, by (2.11)–(2.13) and Markov's inequality, and recalling $M_1/M = A^{-1}n^{-\delta}$,

(2.14)
$$\begin{aligned}\mathbb{P}\left(\sum_{i=1}^M (Y_i - \mathbb{E}Y_i) > \varepsilon M\right) &\leq \varepsilon^{-r} M^{-r} \mathbb{E}\left|\sum_{i=1}^M (Y_i - \mathbb{E}Y_i)\right|^r \\ &\leq C_{10}(\mathbb{E}Y_1^2/M)^{r/2} + C_{11} M^{1-r}\mathbb{E}Y_1^r \\ &\leq C_{12}(M_1/M)^{r/2} + C_{13} M(M_1/M)^r \\ &\leq C_{14} n^{-r\delta/2} + C_{15} n^{1-r\delta} = o(n^{-1}).\end{aligned}$$

For the second sum in (2.10), write $I_i := \mathbf{1}[X_i > M_1]$ and note that for any positive integer $L$, since the $I_i$ are i.i.d.,

(2.15)
$$\begin{aligned}\mathbb{P}\left(\sum_{i=1}^M I_i \geq L\right) &\leq \binom{M}{L}\mathbb{P}(I_1 = \cdots = I_L = 1) \\ &= \binom{M}{L}\mathbb{P}(I_1 = 1)^L \\ &\leq M^L \mathbb{P}(X > M_1)^L.\end{aligned}$$

Moreover, by (2.8) again and the choice of $M$,
$$M\mathbb{P}(X > M_1) \leq C_5 A n^{1/(\gamma-1)+(2-\gamma)(1/(\gamma-1)-\delta)} = CAn^{(\gamma-2)\delta - (\gamma-3)/(\gamma-1)}.$$

Choose $\delta > 0$ such that $\delta_1 := (\gamma-3)/(\gamma-1) - (\gamma-2)\delta > 0$. Thus $M\mathbb{P}(X > M_1) = O(n^{-\delta_1})$ and, by (2.15),

(2.16)
$$\mathbb{P}\left(\sum_{i=1}^M I_i \geq L\right) = O(n^{-L\delta_1}).$$

Now choose $L > 1/\delta_1$, and then $A \geq C_2(L+1)/\varepsilon$. Thus, recalling (1.4),

(2.17) $$\frac{M}{\Delta} \geq \frac{A}{C_2} \geq (L+1)/\varepsilon.$$

It follows from (2.16) that

$$\mathbb{P}\left(\sum_{i=1}^M I_i \geq \varepsilon\frac{M}{\Delta} - 1\right) \leq \mathbb{P}\left(\sum_{i=1}^M I_i \geq L\right) = O(n^{-L\delta_1}) = o(n^{-1}).$$

Consequently, (2.10) shows that, with this choice of $A$, $\mathbb{P}(|\mathcal{C}(v)| > M) = o(n^{-1})$. Since this holds for every vertex $v$, it follows that whp $|\mathcal{C}(v)| \leq M$ for every vertex $v$, and thus $|\mathcal{C}_1| \leq M = An^{1/(\gamma-1)}$.

**3. Proof of Theorem 1.3.** We continue to consider the random multigraph $G^*(n, (d_i)_1^n)$ and prove Theorem 1.3 for this case, which as above implies the result for $G(n, (d_i)_1^n)$ also.

We begin with the upper bound, $|\mathcal{C}_1| \leq (1-\nu)^{-1}\Delta + o_p(n^{1/(\gamma-1)})$. Fix $\varepsilon > 0$. Say that a vertex is *fat* if its degree is at least $M_2 := \varepsilon n^{1/(\gamma-1)}$, and say that a component is *bad* if it contains at least two fat vertices. By (1.3), the number of fat vertices in $G^*(n, (d_i)_1^n)$ is

(3.1) $$n\mathbb{P}(D \geq M_2) \leq C_1 n M_2^{1-\gamma} = C_1 \varepsilon^{1-\gamma} = C_{16}.$$

LEMMA 3.1. *There is whp no bad component.*

PROOF. Let, as in Section 2, $M = An^{1/(\gamma-1)}$, with $A$ so large that whp $|\mathcal{C}_1| \leq M$. Consider a given vertex $v$ and the exploration process in Section 2. By (3.1), the number of fat vertices is at most $C_{16}$ and thus the number of half-edges leading to a fat vertex is at most $C_{16}\Delta$. Consequently, at each step $i$ of the exploration process, the probability of choosing a half-edge $y_i$ leading to a fat vertex $v_i$ is at most $C_{16}\Delta/(n\mu - (2i-1))$ and thus the probability of finding a fat vertex at some step $i \leq M$ is at most, using (1.4),

$$M\frac{C_{16}\Delta}{n\mu - 2M} = O(n^{1/(\gamma-1)+1/(\gamma-1)-1}) = o(1).$$

Thus, for any given $v$, whp the component $\mathcal{C}(v)$ contains no fat vertex found before step $M$, except possibly $v$ itself. On the other hand, we have chosen $M$ such that whp $|\mathcal{C}(v)| \leq M$; hence, for every given $v$, whp $\mathcal{C}(v)$ contains no fat vertex besides possibly $v$.

We apply this starting with a fat vertex $v$. Since the number of fat vertices is bounded by (3.1), it follows that whp none of the components $\mathcal{C}(v)$ with $v$ fat contains a second fat vertex, so none of them is bad. On the other hand, every bad component has to contain a fat vertex $v$, and thus equal $\mathcal{C}(v)$ for some fat $v$. Consequently, whp there is no bad component at all. □

Let $\delta$, $M_1$, $\delta_1$ and $L$ be as in the proof of Theorem 1.1, but replace $M$ by
$$M' := (1-\nu')^{-1}(\Delta + (L+1)\varepsilon n^{1/(\gamma-1)}) \leq C_{17} n^{1/(\gamma-1)}.$$
Further, let $M_2 := \varepsilon n^{1/(\gamma-1)}$ as above. Consider once more the exploration process started at a given vertex $v$. If $|\mathcal{C}(v)| > M'$ and $\mathcal{C}(v)$ is good, then at most one of the degrees $d(v), d(v_1), \ldots, d(v_{M'})$ is greater than $M_2$, and that degree is at most $\Delta$. Consequently, with $I_i := \mathbf{1}[X_i > M_1]$ as above for $i \geq 1$, and $I_0 := 1$ and $v_0 := v$,
$$d(v) + \sum_{i=1}^{M'} I_i \xi_i \leq \sum_{i=0}^{M'} I_i d(v_i) \leq \Delta + \left(\sum_{i=0}^{M'} I_i - 1\right) M_2 = \Delta + M_2 \sum_{i=1}^{M'} I_i.$$
Consequently, in this case (2.9) may be improved to
$$M' \leq \sum_{i=1}^{M'} \xi_i \mathbf{1}[X_i \leq M_1] + \Delta + M_2 \sum_{i=1}^{M'} I_i$$
$$\leq \nu' M' + \sum_{i=1}^{M'} (Y_i - \mathbb{E} Y_i) + \Delta + M_2 \sum_{i=1}^{M'} I_i$$
and thus
$$\sum_{i=1}^{M'} (Y_i - \mathbb{E} Y_i) + M_2 \sum_{i=1}^{M'} I_i \geq (1-\nu') M' - \Delta = (\varepsilon + L\varepsilon) n^{1/(\gamma-1)};$$
hence either
$$\sum_{i=1}^{M'} (Y_i - \mathbb{E} Y_i) \geq \varepsilon n^{1/(\gamma-1)} \geq \varepsilon C_{17}^{-1} M' \tag{3.2}$$
or, recalling $M_2 := \varepsilon n^{1/(\gamma-1)}$,
$$\sum_{i=1}^{M'} I_i \geq L. \tag{3.3}$$
By (2.14) (replacing $\varepsilon$ by $\varepsilon/C_{17}$) and (2.16), both events (3.2) and (3.3) have probabilities $o(n^{-1})$. Consequently,
$$\mathbb{P}(|\mathcal{C}(v)| > M' \text{ and } \mathcal{C}(v) \text{ is good}) = o(n^{-1}).$$
Since this holds for every given vertex $v$, whp $|\mathcal{C}(v)| \leq M'$ for every good component $\mathcal{C}(v)$, which together with Lemma 3.1 shows that whp
$$|\mathcal{C}_1| \leq M' = \frac{\Delta}{1-\nu-\varepsilon} + \frac{L+1}{1-\nu-\varepsilon} \varepsilon n^{1/(\gamma-1)} \leq \frac{\Delta}{1-\nu} + C_{18} \varepsilon n^{1/(\gamma-1)}$$

(with $C_{18}$ not depending on $\varepsilon$). Since $\varepsilon$ can be chosen arbitrarily small, this shows that $|\mathcal{C}_1| \leq (1-\nu)^{-1}\Delta + o_p(n^{1/(\gamma-1)})$, which is the upper bound we sought.

To obtain a corresponding lower bound, let $v^*$ be the vertex with maximum degree (choose any of them if there are several), and consider the component $\mathcal{C}(v^*)$ containing $v^*$.

Let $N_j$ be the number of vertices of distance $j$ from $v^*$. As said in the Introduction, the idea is that (by the usual branching process approximation) $N_j \approx \nu^{j-1} d(v^*) = \nu^{j-1}\Delta$, and summing over $j \geq 1$ yields $\approx (1-\nu)^{-1}\Delta$. In order to make this precise, we for simplicity truncate and consider only generations $j \leq J$ and vertices of degree $\leq K$, where $J$ and $K$ are two fixed large integers. We thus let $N_{jk}$ be the number of vertices of degree $k$ and distance $j$ from $v^*$. We find these by the exploration process in Section 2, now specifying the breadth-first version. Let $\mathcal{H}_j$, $j \geq 0$, be the set of active half-edges when we have processed all half-edges at vertices of distance less than $j$ from $v^*$; thus $\mathcal{H}_j$ is the set of half-edges that may lead to a vertex of distance $j+1$. We have $|\mathcal{H}_0| = d(v^*) = \Delta$.

Let, for $k \geq 1$,
$$p_k^* := \frac{k\mathbb{P}(D=k)}{\mu} = \frac{kn_k}{\mu n},$$

the fraction of half-edges that are adjacent to a vertex of degree $k$ (this is the size-biased distribution of $D$), and note that, by (1.1) and (1.2),
$$\sum_{k=1}^{\infty} p_k^* = 1, \qquad \sum_{k=1}^{\infty}(k-1)p_k^* = \nu.$$

Assume $\nu > 0$ and let $\varepsilon > 0$ be so small that $\varepsilon < \nu$. (The case $\nu = 0$ is simpler and is left to the reader.) Let $\nu_K := \sum_1^K (k-1)p_k^*$ and assume that $K$ is so large that $\nu_K > \varepsilon$ and $\sum_1^K p_k^* > 1 - \varepsilon$. Further assume that $\Delta \geq (\log n)^2$ (otherwise the result is trivial), and assume below that $n$ is large whenever needed.

We may assume that the exploration process stops before it reaches $An^{1/(\gamma-1)}$ steps for some large $A$, either by Theorem 1.1 or because the lower bound otherwise is trivial. Then, at every step $i$ there are $n_k - O(n^{1/(\gamma-1)}) = n\mathbb{P}(D=k) - o(n)$ unused vertices of degree $k$ and thus the probability to connect to one of them equals
$$\frac{kn\mathbb{P}(D=k) - o(n)}{\mu n - o(n)} = p_k^* + o(1) \geq p^* - \varepsilon/K^2,$$

given any history of the exploration. Hence, for every fixed $j$ and $k$ and given $|\mathcal{H}_{j-1}| = h_{j-1}$, $N_{jk}$ dominates a sum of independent indicators with the binomial distribution $\text{Bi}(h_{j-1}, p_k^* - \varepsilon/K^2)$, and thus by the law of large numbers, whp

$$N_{jk} \geq h_{j-1}(p_k^* - 2\varepsilon/K^2), \tag{3.4}$$

provided $h_{j-1} \geq \log n$, say. Moreover, the probability that a given vertex of degree $k$ is chosen twice in the $j$th generation is $O(kh_{j-1}/n)^2 = O(n^{1/(\gamma-1)}/n)^2 = o(n^{-1})$, so whp there are no repetitions and all these $N_{jk}$ vertices yield $k-1$ new active half-edges each.

We consider only a finite number of $(j,k)$, namely $j \leq J$ and $k \leq K$, and thus whp the statements above hold for all these $(j,k)$ simultaneously. Consequently, whp, for all $j \leq J$,

$$|\mathcal{H}_j| \geq \sum_{k=1}^{K}(k-1)N_{jk} \geq |\mathcal{H}_{j-1}|\sum_{k=1}^{K}(k-1)(p_k^* - 2\varepsilon/K^2) \geq |\mathcal{H}_{j-1}|(\nu_K - \varepsilon),$$

and thus by induction (which also verifies $h_{j-1} \geq \log n$)

$$|\mathcal{H}_j| \geq |\mathcal{H}_0|(\nu_K - \varepsilon)^j = (\nu_K - \varepsilon)^j \Delta, \qquad 0 \leq j \leq J.$$

By (3.4) again, it follows that whp, for $1 \leq j \leq J$ and $1 \leq k \leq K$,

$$N_{jk} \geq (\nu_K - \varepsilon)^{j-1}\Delta(p_k^* - 2\varepsilon/K^2),$$

and thus, summing over these $j$ and $k$,

$$|\mathcal{C}(v^*)| \geq \sum_{j=1}^{J}(\nu_K - \varepsilon)^{j-1}\sum_{k=1}^{K}(p_k^* - 2\varepsilon/K^2)\Delta \geq \frac{1 - (\nu_K - \varepsilon)^J}{1 - \nu_K + \varepsilon}(1 - 3\varepsilon)\Delta.$$

For any small $\eta > 0$, we may choose $\varepsilon$ small and $J$ and $K$ large so that this yields, recalling (1.4),

$$|\mathcal{C}_1| \geq |\mathcal{C}(v^*)| \geq \left(\frac{1}{1-\nu} - \eta\right)\Delta \geq \frac{\Delta}{1-\nu} - \eta C_2 n^{1/(\gamma-1)}$$

whp, which yields the required lower bound and completes the proof.

**4. Extensions to other random graphs.** Theorems 1.1 and 1.3 can easily be extended to other random graph models. We will in this section study some models where we can obtain results as corollaries of the results above for $G(n, (d_i)_1^n)$. (Another approach would be to try to adapt the proofs above to other models, but we will not pursue this possibility here.) Indeed, suppose that $G_n$ is a random graph with $n$ vertices (labeled $1, \ldots, n$), such that $G_n$ conditioned on the degree sequence $(d_i(G_n))_{i=1}^n$ is uniformly distributed over all graphs with this degree sequence. Equivalently, $\mathbb{P}(G_n = H_1) = \mathbb{P}(G_n = H_2)$ for any two graphs $H_1$ and $H_2$ with vertices $1, \ldots, n$ and the same degree sequence. (In statistical terminology, the degree sequence is a sufficient statistic.) We then can apply the theorems in Section 1 to $G_n$ conditioned on the degree sequence.

Define, in analogy with the notation above,

$$\mu(G_n) = \frac{1}{n} \sum_{i=1}^{n} d_i(G_n), \tag{4.1}$$

$$\nu(G_n) = \frac{1}{n\mu(G_n)} \sum_{i=1}^{n} d_i(G_n)(d_i(G_n) - 1), \tag{4.2}$$

$$\Delta(G_n) := \max_i d_i(G_n), \tag{4.3}$$

and note that these are random variables depending on the degree sequence of $G_n$ only. Further, let $\mathcal{C}_1(G_n)$ be the largest component of $G_n$ (with any choice among them if there is a tie).

THEOREM 4.1. *Let $\gamma > 3$. Suppose that, for $n \geq 1$, $G_n$ is a random graph with $n$ vertices such that $G_n$ conditioned on the degree sequence is uniformly distributed over all graphs with this degree sequence. Suppose further that for every $\varepsilon > 0$ there exists $\delta > 0$ such that*

$$\mathbb{P}\left(\sup_{k \geq 1} k^{\gamma-1} |\{i : d_i(G_n) \geq k\}|/n > \delta^{-1}\right) < \varepsilon, \tag{4.4}$$

$$\mathbb{P}(\mu(G_n) < \delta) < \varepsilon, \tag{4.5}$$

$$\mathbb{P}(\nu(G_n) > 1 - \delta) < \varepsilon. \tag{4.6}$$

*(In other words, $\sup_{k \geq 1} k^{\gamma-1} |\{i : d_i(G_n) \geq k\}|/n$, $1/\mu(G_n)$ and $1/(1-\nu(G_n))_+$ are stochastically bounded.) Then*

$$|\mathcal{C}_1(G_n)| = \frac{\Delta(G_n)}{1 - \nu(G_n)} + o_p(n^{1/(\gamma-1)}), \tag{4.7}$$

*and, in particular, $|\mathcal{C}_1(G_n)| = O_p(n^{1/(\gamma-1)})$.*

PROOF. First, by conditioning on the degree sequence, we may assume that the events

$$\sup_{k \geq 1} k^{\gamma-1} |\{i : d_i(G_n) \geq k\}|/n \leq \delta^{-1}, \tag{4.8}$$

$$\mu(G_n) \geq \delta, \tag{4.9}$$

$$\nu(G_n) \leq 1 - \delta \tag{4.10}$$

hold (surely) for some $\delta > 0$. Indeed, if we have proved this case, we may condition on these events and see that the conclusion then holds; by (4.4)–(4.6), this leaves only an event of probability $< 3\varepsilon$, which can be made arbitrarily small, and the conclusion follows for $G_n$ without conditioning also.

Second, the Cauchy–Schwarz inequality implies, using (4.10),

$$(\mu(G_n))^2 \leq \frac{1}{n}\sum_{i=1}^{n} d_i(G_n)^2 = \mu(G_n)\nu(G_n) + \mu(G_n) \leq 2\mu(G_n),$$

and thus $\mu(G_n) \leq 2$. Hence, both $\mu(G_n)$ and $\nu(G_n)$ are sequences of random variables that are uniformly bounded, and thus tight; thus there exist subsequences such that $\mu(G_n)$ and $\nu(G_n)$ converge in distribution. We now use the standard subsubsequence trick: if every subsequence of $(G_n)$ has a subsubsequence for which (4.7) holds, then (4.7) holds for the full sequence. Consequently, by considering subsequences, we may assume that $\mu(G_n) \xrightarrow{d} \mu_\infty$ and $\nu(G_n) \xrightarrow{d} \nu_\infty$ for some random variables $\mu_\infty$ and $\nu_\infty$, with $\delta \leq \mu_\infty \leq 2$ and $0 \leq \nu_\infty \leq 1-\delta$.

Third, by the Skorohod coupling theorem ([13], Theorem 4.30), we may assume that all random graphs $G_n$ are defined on a common probability space and that the limits $\mu(G_n) \to \mu_\infty$ and $\nu(G_n) \to \nu_\infty$ hold a.s. and not just in distribution. We now condition $G_n$ on its degree sequence $(d_i(G_n))_{i=1}^n$. By the reductions above, the assumptions of Theorem 1.1 are satisfied a.s. Hence Theorem 1.3 shows that for every $\eta > 0$

$$(4.11) \quad \mathbb{P}\bigg(\bigg||\mathcal{C}_1(G_n)| - \frac{\Delta(G_n)}{1-\nu(G_n)}\bigg| > \eta n^{1/(\gamma-1)} \bigg| (d_i(G_n))_{i=1}^n\bigg) \to 0 \qquad \text{a.s.},$$

and the result follows by taking the expectation in (4.11) (using dominated convergence). $\square$

As an illustration, we apply this theorem to two popular models of random graphs.

4.1. *The Newman, Strogatz, Watts model.* Durrett [7] actually stated his conjecture for the following random graph model, introduced by Newman, Strogatz and Watts [17, 18]. Let $D$ be a nonnegative integer-valued random variable. Take i.i.d. copies $D_1, D_2, \ldots$ of $D$. Condition on $\sum_{i=1}^n D_i$ being even, and let $G_n$ be the random multigraph $G^*(n,(D_i)_1^n)$, defined by the configuration model as in Section 2. (This is a multigraph. We may either accept this, or delete all loops and merge parallel edges, or condition on the graph being simple; this does not affect our result.)

COROLLARY 4.2. *Let $\gamma > 3$. Assume that $\mathbb{P}(D \geq k) = O(k^{1-\gamma})$ and $\nu := \mathbb{E}D(D-1)/\mathbb{E}D < 1$. Then*

$$|\mathcal{C}_1(G_n)| = \frac{\Delta(G_n)}{1-\nu} + o_p(n^{1/(\gamma-1)}).$$

REMARK 4.3. If further, as assumed by Durrett [7], $\mathbb{P}(D=k) \sim ak^{-\gamma}$ as $k \to \infty$, with $a > 0$, then $\mathbb{P}(D \geq k) \sim a_1 k^{1-\gamma}$ with $a_1 = a/(\gamma - 1)$, and it is easily seen that $n^{-1/(\gamma-1)} \Delta(G_n) \xrightarrow{d} a_1^{1/(\gamma-1)} Z$, with $\mathbb{P}(Z \leq x) = \exp(-x^{1-\gamma})$ for $x > 0$, one of the classical extreme value distributions. (See, e.g., Leadbetter, Lindgren and Rootzén [14], Chapter 1 and Example 1.7.6.) Consequently,

$$n^{-1/(\gamma-1)} |\mathcal{C}_1(G_n)| \xrightarrow{d} \frac{a_1^{1/(\gamma-1)}}{1-\nu} Z.$$

PROOF OF COROLLARY 4.2. Let $\mu = \mathbb{E}D$. The case $\mu = 0$ is trivial, so we may assume $\mu > 0$. The assumption $\mathbb{P}(D \geq k) = O(k^{1-\gamma})$ with $\gamma > 3$ implies that $\mathbb{E}D^2 < \infty$, so $0 < \mu < \infty$ and $0 \leq \nu < \infty$. The (weak) law of large numbers shows that $\sum_{i=1}^{n} D_i/n \xrightarrow{P} \mu$ and $\sum_{i=1}^{n} D_i(D_i - 1)/n \xrightarrow{P} \mathbb{E}D(D-1) = \mu\nu$, and these hold also if we condition on $\sum_{i=1}^{n} D_i$ even, which has probability converging to $1/2$ [or 1, if $\mathbb{P}(D_i \text{ odd}) = 0$]. Hence, using the notation (4.1) and (4.2), $\mu(G_n) \xrightarrow{P} \mu$ and $\nu(G_n) \xrightarrow{P} \nu$, which yields (4.5) and (4.6) for any $\delta < \min(\mu, 1 - \nu)$. This also implies

$$\frac{\Delta(G_n)}{1-\nu(G_n)} - \frac{\Delta(G_n)}{1-\nu} = \frac{\nu(G_n) - \nu}{(1-\nu)(1-\nu(G_n))} \Delta(G_n) = o_p(n^{1/(\gamma-1)}),$$

so we may replace $\nu(G_n)$ by $\nu$ in (4.7).

It remains only to verify (4.4). Again, it suffices to show the corresponding result for the i.i.d. sequence $(D_i)_1^n$ before conditioning. Let $N_n(x) := |\{i \leq n : D_i \geq x\}|$ and let $\mathcal{E}(B, x)$ be the "bad" event $\{N_n(x) > Bnx^{1-\gamma}\}$, for real $x, B > 0$. Then $N_n(x) \sim \text{Bi}(n, \mathbb{P}(D > x))$, and by assumption $\mathbb{P}(D > x) \leq C_{19} x^{1-\gamma}$. Assume that $B \geq 2C_{19}$, so $\mathbb{E}N_n(x) = n\mathbb{P}(D > x) \leq \frac{1}{2}Bnx^{1-\gamma}$ and $\text{Var } N_n(x) \leq \mathbb{E}N_n(x) \leq \frac{1}{2}Bnx^{1-\gamma}$. Then, by Chebyshev's inequality,

$$(4.12) \quad \mathbb{P}(\mathcal{E}(B,x)) \leq \frac{\text{Var}(N_n(x))}{(Bnx^{1-\gamma} - \mathbb{E}N_n(x))^2} \leq \frac{\text{Var}(N_n(x))}{((1/2)Bnx^{1-\gamma})^2} \leq \frac{2}{Bnx^{1-\gamma}}.$$

Also, since $N_n(x)$ is integer-valued,

$$(4.13) \quad \mathbb{P}(\mathcal{E}(B,x)) \leq \mathbb{P}(N_n(x) \geq 1) \leq \mathbb{E}(N_n(x)) \leq C_{19} n x^{1-\gamma}.$$

We consider a geometric sequence and let $x_j := (2^j n)^{1/(\gamma-1)}$, $j \in \mathbb{Z}$. By (4.12) and (4.13),

$$(4.14) \quad \mathbb{P}(\mathcal{E}(B, x_j)) \leq \max\left(\frac{2^{j+1}}{B}, C_{19} 2^{-j}\right).$$

Hence, each $\mathbb{P}(\mathcal{E}(B, x_j)) \to 0$ as $B \to \infty$; moreover, (4.14) further shows that $\mathbb{P}(\mathcal{E}(B, x_j)) \leq C_{20} 2^{-|j|}$. Consequently, dominated convergence shows that

$$\sum_{j=-\infty}^{\infty} \mathbb{P}(\mathcal{E}(B, x_j)) \to 0 \quad \text{as } B \to \infty.$$

In other words, by choosing $B$ large enough, with probability at least $1-\varepsilon$ none of $\mathcal{E}(B, x_j)$ holds, that is, $N_n(x_j) \leq Bnx_j^{1-\gamma}$ for all $j \in \mathbb{Z}$, which implies $N_n(x) \leq 2Bnx^{1-\gamma}$ for all $x > 0$. [Use $N_n(x) \leq N_n(x_j)$ for the largest $x_j \leq x$.] This verifies (4.4), and the result follows by Theorem 4.1. $\square$

A minor variation of this model has been used by van der Hofstad, Hooghiemstra and Van Mieghem [22]; the same result holds for their version.

4.2. *An inhomogeneous Erdős–Rényi graph.* A model that has been used, with minor variations, by several authors is the following: Let $W_1, \ldots, W_n$ be nonnegative real numbers, random or deterministic. Given these numbers, let

$$p_{ij} := \frac{W_i W_j}{n + W_i W_j}, \tag{4.15}$$

and let $G_n$ be the random graph on $n$ vertices where [conditioned on $(W_i)_1^n$] edges appear independently, with $\mathbb{P}(\text{an edge } ij) = p_{ij}$ for $1 \leq i < j \leq n$. It was observed by Britton, Deijfen and Martin-Löf [3] that with the definition (4.15), this yields a random graph that conditioned on its degree sequence is uniformly distributed, as we require here.

A version of this model, which may look simpler, is to replace (4.15) by $p_{ij} := W_i W_j / n$ [or rather, when necessary, $\min(W_i W_j / n, 1)$]. However, this ruins the uniformity after conditioning on the degree sequence. Nevertheless, the difference between the two versions is small, unless some $W_i$ are very large, and usually the difference is negligible. More precisely, it can be shown [10] that if $\sum_{i=1}^n W_i^3 = o_p(n^{3/2})$, which for example is the case if, as below, $W_i$ are i.i.d. with $\mathbb{E}W_i^2 < \infty$, then the two versions are asymptotically equivalent in a strong sense (the total variation distance tends to 0), and thus the result below holds for this version also.

Another version, introduced by Chung and Lu [4, 5], is to take $p_{ij} := W_i W_j / \sum_{k=1}^n W_k$. This becomes the same as the preceding version if we replace $W_i$ by $W_i (\sum_k W_k / n)^{-1/2}$. Again, results are easily transferred to this version; we omit the details. Another minor (and for us inessential) variation is used by Norros and Reittu [19], and several versions are studied by van den Esker, van der Hofstad and Hooghiemstra [21].

It is common to take $(W_i)_1^n$ i.i.d. (e.g., this is used in [3, 19, 21]), but it is also possible to use suitable deterministic sequences [4, 5] or dependent random variables; see Bollobás, Janson and Riordan [2], Section 16.4, for a general discussion and further references. (The model is a special case of the inhomogeneous random graph defined in [2].) For simplicity we treat here only the i.i.d. case. The condition for existence of a giant component then simply is $\mathbb{E}W_i^2 > 1$; see Bollobás, Janson and Riordan [2], Section 16.4.

COROLLARY 4.4. *Let $G_n$ be the random graph defined as above, using (4.15) where $W_i$ are i.i.d. copies of a given nonnegative random variable $W$. Let $\gamma > 3$. Assume that $\mathbb{P}(W \geq x) = O(x^{1-\gamma})$, $x > 0$, and $\nu := \mathbb{E}W^2 < 1$. Then*

$$|\mathcal{C}_1(G_n)| = \frac{\Delta(G_n)}{1-\nu} + o_p(n^{1/(\gamma-1)}).$$

PROOF. We may assume $\mathbb{E}W > 0$. Since $n\mu(G_n) = \sum_i d_i(G_n)$ is twice the number of edges in $G_n$, it is easily seen that $\mu(G_n) \xrightarrow{\mathrm{P}} \mu := (\mathbb{E}W)^2 > 0$ which verifies (4.5); for example, this follows from the more general Lemma 8.1 in [2].

Similarly, for example by observing that $\sum_i d_i(G_n)(d_i(G_n) - 1)$ is twice the number of paths of length 2 in $G_n$ and using [2], Theorem 17.1,

$$\mu_n(G_n)\nu_n(G_n) = \frac{1}{n}\sum_{i=1}^n d_i(G_n)(d_i(G_n) - 1) \xrightarrow{\mathrm{P}} (\mathbb{E}W)^2 \mathbb{E}W^2.$$

Hence, $\nu(G_n) \xrightarrow{\mathrm{P}} \nu := \mathbb{E}W^2$, which verifies (4.6). This also enables us to replace $\nu(G_n)$ by $\nu$ in (4.7).

It remains to verify (4.4). First, by the argument in the proof of Corollary 4.2, for any $\varepsilon > 0$,

(4.16)  $\quad |\{i \leq n : W_i \geq k\}| \leq Bnk^{1-\gamma} \quad$ for all $k \geq 1$,

with probability $> 1 - \varepsilon$, provided $B = B(\varepsilon)$ is large enough. Hence, conditioning on (4.16), we may assume that (4.16) holds. We condition on $(W_i)_1^n$ and assume thus in the remainder of the proof that (for each $n$) $W_i$ are fixed numbers such that (4.16) holds, for some constant $B$. In particular, this implies

(4.17) $\quad \sum_{i=1}^n W_i = \sum_{k=1}^\infty |\{i \leq n : W_i \geq k\}| \leq \sum_{k=1}^\infty Bnk^{1-\gamma} \leq C_{21} n.$

Each $d_i(G_n)$ is a sum of independent 0–1 variables, and by (4.15) and (4.17)

$$\mathbb{E}d_i(G_n) = \sum_{j \neq i} \frac{W_i W_j}{n + W_i W_j} \leq \frac{W_i \sum_j W_j}{n} \leq C_{21} W_i.$$

As a consequence, we have the Chernoff estimate (see [12], Corollary 2.4 and Theorem 2.8)

(4.18) $\quad \mathbb{P}(d_i(G_n) \geq t) \leq \exp(-t), \qquad t \geq 7 C_{21} W_i.$

Let $j \geq 1$ and let

$$I_j := \{i \leq n : W_i \leq (7C_{21})^{-1} 2^j\}, \qquad I'_j := \{1,\ldots,n\} \setminus I_j.$$

By (4.16), $|I'_j| \leq C_{22} 2^{(1-\gamma)j} n$. Further, let $I^*_j$ be the random set $\{i \in I_j : d_i(G_n) \geq 2^j\}$. By (4.18),

$$\mathbb{E}|I^*_j| \leq \exp(-2^j)|I_j| \leq \exp(-2^j)n \leq C_{23} 2^{-\gamma j} n,$$

and thus, for every $\varepsilon > 0$,

$$\mathbb{P}(|I^*_j| > \varepsilon^{-1} C_{23} 2^{(1-\gamma)j} n) \leq \varepsilon 2^{-j}.$$

Hence, with probability $\geq 1 - \varepsilon$, $|I^*_j| \leq \varepsilon^{-1} C_{23} 2^{(1-\gamma)j} n$ for all $j \geq 1$, and then, assuming $\varepsilon \leq 1$,

$$|\{i \leq n : d_i(G_n) \geq 2^j\}| \leq |I^*_j| + |I'_j| \leq \varepsilon^{-1} C_{24} 2^{(1-\gamma)j} n, \qquad j \geq 1,$$

which implies (4.1) and completes the proof. $\square$

**Acknowledgments.** This research was done mainly during the 11th Brazilian School of Probability in Maresias, August 2007. I thank Rick Durrett for interesting discussions and comments.

DEPARTMENT OF MATHEMATICS
UPPSALA UNIVERSITY
PO BOX 480
SE-751 06 UPPSALA
SWEDEN
E-MAIL: svante.janson@math.uu.se
URL: http://www.math.uu.se/~svante/